\documentclass [a4paper, 12t]{article}
\usepackage{amsmath}
\usepackage{amsfonts}
\usepackage{amssymb}
\usepackage{amsthm}
\usepackage{makeidx}
\usepackage{graphicx}
\usepackage{pb-diagram}
\usepackage{epstopdf}
\usepackage{fancyhdr}
\usepackage{enumitem}
\usepackage[all]{xy}

\newtheorem{teo}{Theorem}
\newtheorem{prop}[teo]{Proposition}
\newtheorem{lem}{Lemma}

\theoremstyle{definition}

\newtheorem{con*}{Conjecture}

\begin{document}
\setlength{\parskip}{1ex plus 0.5ex minus 0.2ex}
\begin{center}

\textbf{A COUNTEREXAMPLE FOR A SUP THEOREM IN LOCALLY CONVEX SPACES}\\
\end{center}
\begin{center}
\textsl{By STEFANO ROSSI}
\end{center}
$$$$
$$$$
\textbf{Abstract}\quad In this brief note, we provide an example of non complete locally convex space $E$ with a $\sigma(E,E^*)$ closed bounded subset $\mathcal{C}\subset E$, which is not $\sigma(E, E^*)$-compact, even if every $\varphi\in E^*$ attains its sup over $\mathcal{C}$.
$$$$
$$$$
In \cite{Jam2} James has shown the following far reaching theorem, that can be viewed as one the deepest result on weak topologies together with the Eberlein-Smulian theorem:
\begin{teo}
Let $E$ be a \emph{complete} locally convex space. If $\mathcal{C}\subset E$ is a bounded $\sigma(E, E^*)$-closed subset, then the following are equivalent:
\begin {enumerate}
\item $\mathcal{C}$ is $\sigma(E,E^*)$-compact.
\item Given any $\varphi\in E^*$, there is $x\in\mathcal{C}$ such that $\sup_{y\in\mathcal{C}} |\varphi(y)|=|\varphi(x)|$.
\end{enumerate}
\end{teo}
Completeness assumption cannot be dropped as R.C. James himself has pointed out in \cite{Jam1}, where an example of non complete normed space whose functionals are all norm-attaining is given. Clearly, by James' theorem the completion of this space is a (not stricly convex) reflexive Banach space.\\
As mentioned in the abstract, in this note we will give an example of a non complete locally convex space $E$ with a convex bounded $\sigma(E,E^*)$-closed subset $\mathcal{C}\subset E$, which is not weak compact, though every $\varphi\in E^*$ has a maximum on $\mathcal{C}$.\\

Before showing our counterexample, we start recalling some basic definitions.
Given a Banach space $\mathfrak{X}$, a bounded linear functional $\varphi\in\mathfrak{X}^*$ is said to be \emph{norm-attaining} if there exists $x\in\mathfrak{X}$ with $\|x\|=1$ such that $\|\varphi\|=|\varphi(x)|$. Obviously, a subspace $\mathfrak{M}\subset\mathfrak{X}^*$ is norm-attaining if each functional $\varphi\in\mathfrak{M}$ is norm-attaining. Finally a subspace $\mathfrak{M}\subset\mathfrak{X}^*$ is \emph{determinant} (or $1$-norming) if $\|x\|= sup_{\varphi\in\mathfrak{M}_1} |\varphi(x)|$ for each $x\in\mathfrak{X}$, $\mathfrak{M}_1$ being the unit ball of $\mathfrak{M}$.\\

Let $\mathfrak{X}$ be the (complex) Banach space of those functions $f:[0,1]\subset\mathbb{R}\rightarrow\mathbb{C}$ such that $\sum_{x\in [0,1]}|f(x)|<\infty$, endowed with the norm $\|\cdot\|_1$ given by $\|f\|_1=\sum_{x\in [0,1]}|f(x)|$. Note that the support of any $f\in\mathfrak{X}$ is at most countable.\\
Let $f_x\in\mathfrak{X}$ be the function given by $f_x(y)=\delta_{x,y}$ for each $y\in[0,1]$. Since $$\|f_x-f_y\|_1=2\quad\textrm{for each}\, x,y\in [0,1]$$ $\mathfrak{X}$ cannot be separable, $[0,1]$ being uncontable.\\
$B[0,1]$ will indicate the space of all real bounded function defined on $[0,1]$, endowed with the sup-norm $\|\cdot\|_{\infty}$. The following simple lemma recognizes $B[0,1]$ as the dual space of $\mathfrak{X}$.
\begin{lem}
$\Phi:B[0,1]\rightarrow \mathfrak{X}^*$, given by $\Phi(g)(f)=\sum_{x\in [0,1]}g(x)f(x)$ for each $g\in B[0,1]$ and $f\in\mathfrak{X}$, is an isometric isomorphism.
\end{lem}

\begin{proof}
Clearly $\Phi(g)(f)=\sum_{x\in [0,1]}g(x)f(x)$ defines a bounded linear map from $B[0,1]$ to $\mathfrak{X}^*$ and $\|\Phi(g)\|\leq\|g\|_{\infty}$ for each $g\in B[0,1]$.\\
Now,  let $\varphi\in\mathfrak{X}^*$. Let us define $g$ as the function given by $g(x)=\varphi( f_x)$ for each $x\in [0,1]$. Since $|g(x)|\leq\|\varphi\|$, we get $g$ is a bounded function with $\|g\|_{\infty}\leq\|\varphi\|$. To conclude the proof, it only remains to check that $\varphi=\Phi(g)$. If $f\in\mathfrak{X}$, we can write $f=\sum_{i=1}^{\infty}f(x_i)f_{x_i}$, where $\{x_i:i\in\mathbb{N}\}$ is the support of $f$. Put $f_n=\sum_{i=1}^{n}f(x_i)f_{x_i}$, clearly we have $\|f-f_n\|_1\rightarrow 0$, hence: $$\varphi (f)=\lim_{n\to\infty}\varphi(f_n)=\lim_{n\to\infty}\sum_{i=1}^{n}f(x_i)g(x_i)=\sum_{i=1}^{\infty}f(x_i)g(x_i)=\Phi(g)(f)$$
This concludes the proof, since $f\in\mathfrak{X}$ is arbitrary.
\end {proof}
Now let us consider the space $C[0,1]$ of continuous functions on $[0,1]$ as a closed subspace of $\mathfrak{X}^*=B[0,1]$.
We have the following:
\begin{prop}
$C[0,1]\subset\mathfrak{X}^*$ is a determinant subspace of norm-attaining linear forms.
\end{prop}
\begin{proof}
Let $f\in\mathfrak{X}$ and let $\{x_i:i\in\mathbb{N}\}$ be its support. Given any $\varepsilon>0$ there is a natural number $N$ such that $\sum_{N+1}^{\infty}|f(x_i)|<\varepsilon$. Now let $g$ be any continuous function on $[0,1]$ such that $g(x_i)=e^{i\theta_i}$\, for each $i=1,2\dots,N$ and $\|g\|_{\infty}=1$, where $\theta_i\in\mathbb{R}$ is given by $f(x_i)=|f(x_i)|e^{-i\theta_i}$. We have:
$$\Phi(g)(f)=\sum_{i=1}^N |f(x_i)|+\sum_{i=N+1}^{\infty} f(x_i)g(x_i)> \|f\|-2\varepsilon$$
since $|\sum_{i=N+1}^{\infty} f(x_i)g(x_i)|<\varepsilon$. This proves that $$\|f\|_1=\sup_{g\in C[0,1]: \|g\|_{\infty}=1} |\Phi(g)(f)|$$
so that $C[0,1]$ is a determinant subspace.\\
Now, let $g\in\ C[0,1]$, we have $\|\Phi(g)\|=\|g\|_{\infty}=|g(x_0)|$ for some $x_0\in [0,1]$. Then $\|\Phi(g)\|=\Phi(g)(\tilde{f})$, where $\tilde{f}=e^{i\theta_0}f_{x_0}$ ($g(x_0)=e^{-i\theta_0}|g(x_0)|$). This concludes the proof.
\end{proof}
Now we consider the dual pair (cfr. \cite{Wil}) with the natural pairing (evaluation). $(\mathfrak{X},\mathfrak{M})$ and the locally convex space $E$, which is $\mathfrak{X}$ endowed with the Mackey topology $\tau(\mathfrak{X},\mathfrak{M})$, so that we have $E^*=\mathfrak{M}$ by virtue of the Mackey-Arens theorem (see \cite{Wil}). In \cite{Josè} it is shown that $E$ is nor complete, neither quasi-complete.\\
Let $\mathcal{C}\subset E$ be the (circled) convex subset given by the unit ball of $\mathfrak{X}_1$. Since $\mathfrak{M}$ is determinant, we have: $$\mathcal{C}=\bigcap_{\varphi\in\mathfrak{M}_1\subset E^*} \{f\in E: |\varphi(f)|\leq 1\}$$
thus $\mathcal{C}$ is $\sigma(E,E^*)$-closed.
Clearly $\mathcal{C}$ is a Mackey bounded subset as a norm-bounded subset, anyway it is not $\sigma(E, E^*)$- compact, even if each $\varphi\in E^*$ has a maximum on $\mathcal{C}$:
\begin{prop}
$\mathcal{C}\subset E$ is not $\sigma(E, E^*)$-compact.
\end{prop}
\begin{proof}
We argue by \emph{reductio ad absurdum}. Let us consider the linear map $F: \mathfrak{X}\rightarrow C[0,1]^*$ given by $F(f)(g)= \sum_{x\in [0,1]} f(x)g(x)$ for each $g\in C[0,1]$, for each $f\in\mathfrak{X}$. Since $C[0,1]$ is determinant (for $\mathfrak{X}$), $F$ is an isometry; moreover it is a $\sigma(\mathfrak{X},\mathfrak{M})$-$\sigma(\mathfrak{M}^*,\mathfrak{M})$ continuous map ($\mathfrak{M}^*$ is the Banach dual space of $\mathfrak{M}$).\\
Observe that $Ran F$ is weak*-dense in $C[0,1]^*$, since one trivially has $Ran F^{\perp}=0$.
If $\mathcal{C}=\mathfrak{X}_1$ is $\sigma(E, E^*)$ compact, then $F(C)=F(\mathfrak{X})_1$ is weakly* closed, so $Ran F$ is weakly* closed thanks to Krein Smulian theorem. This means that $F$ is surjective, that is $C[0,1]^*\cong\mathfrak{X}$. This is absurd, since $C[0,1]^*=\mathcal{M}([0,1])$ (finite Borel measures on $[0,1]$) and $\mathfrak{X}\subset\mathcal{M}([0,1])$ as a \emph{proper} sunbpace, where the (isometric) inclusion map is simply given by $i(f)=\sum_{x\in [0,1]}f(x)\delta_x$, $\delta_x$ being the Dirac measure concetrated on $x$.
\end{proof}

$$$$
$$$$
\emph{DIP. MAT. GUIDO CASTELNUOVO UNIV. DI ROMA LA SAPIENZA, ROME, ITALY}\\
\emph{E-mail address}: \verb"s-rossi@mat.uniroma1.it"

\end{document}